\documentclass[a4paper,11pt]{article}
\pdfoutput=1 

\usepackage{jheppub} 

\usepackage[T1]{fontenc} 

\title{\boldmath Regions exhibiting stationary values of average cross section}

\author[a]{Kyeong Min Kim}

\affiliation[a]{Mathematical Institute, University of Oxford, Woodstock Road, Oxford, OX2 6GG, UK}

\emailAdd{kyeong.kim@st-hughs.ox.ac.uk}

\abstract{The problem of maximizing the average cross section through a point within a shape is introduced. This idea is extended into arbitrary dimensions. However, the average cross sectional volume cannot be maximized unless the cross sections pass through the centroid of the shape. Therefore, we focus on the shapes with stationary values of average cross section. A novel formula for average cross section in any dimension is presented. The equation of the stationary shape is derived and the general characteristics are discussed through algebraic solution, numerical analysis, and interpretation with graphical display. The special properties of the stationary shapes in 2, 3, and 5 dimension is examined.}

\begin{document} 
\maketitle
\flushbottom

\section{Introduction}

The isoperimetric problem determines that the area of the two-dimensional shape with its boundary as fixed length is maximized when the shape is the disc. \cite{A,B} Similarly, the ball is known to be the shape with its surface area minimized when its volume is fixed. \cite{C} The same would be true for a sphere in any dimensions. The sphere in $n$-dimensions, $S^n$, is defined as the set of points in a distance $r$ from a central point, or
\begin{equation}
    S^n=\{x\in R^{n+1} : \lvert x\rvert = r \}.
\end{equation}
These properties of the spheres, in whatever dimensions, can be obtained through the calculus of variation.
These properties of the spheres, in whatever dimensions, can be obtained through the calculus of variation.

The initial idea of the problem that this paper explores is described in elementary geometry analogous to the isoperimetric problem. It is a natural idea of how to maximize the average cross section of the shape through the centroid. Consider the case where we maximize the average diameter in two-dimensional shapes. The only constraint of this problem is that the area of the shape is fixed. We do not constrain the shape to be star shaped or connected.

The solution to this problem is the disc, and it is not a surprising result, because spheres exhibit unique features. The exact reason why the disc gives the maximum average diameter is given in the Section 4, after the derivation of the average cross sectional thickness.

The problem can be extended in two aspects. We can consider the case when the cross sections of a shape pass through any point inside the shape rather than through centroid. Also, the problem can be broadened to arbitrary dimensions of shapes and cross sections.

Firstly, the conditions on the cross sectional areas can be altered. Let the cross sections pass through a point which is at a certain distance away from the centre of mass. Set the coordinate system, such that the cross sections pass through the origin of the coordinate. Then the constraint that the centre of mass of the two-dimensional shape is
\begin{equation}
    \frac{\underline{M}}{A}=\frac{1}{A}\int \underline{x} \,dA
\end{equation}
away from the origin is added to the initial problem, where $\underline{M}$ is the moment and $A$ is the area of the shape. In this extended version of the problem, the sphere would not be the shape that maximizes its average cross sectional area.

\begin{figure}[tbp]
\centering 
\includegraphics[width=.8\textwidth]{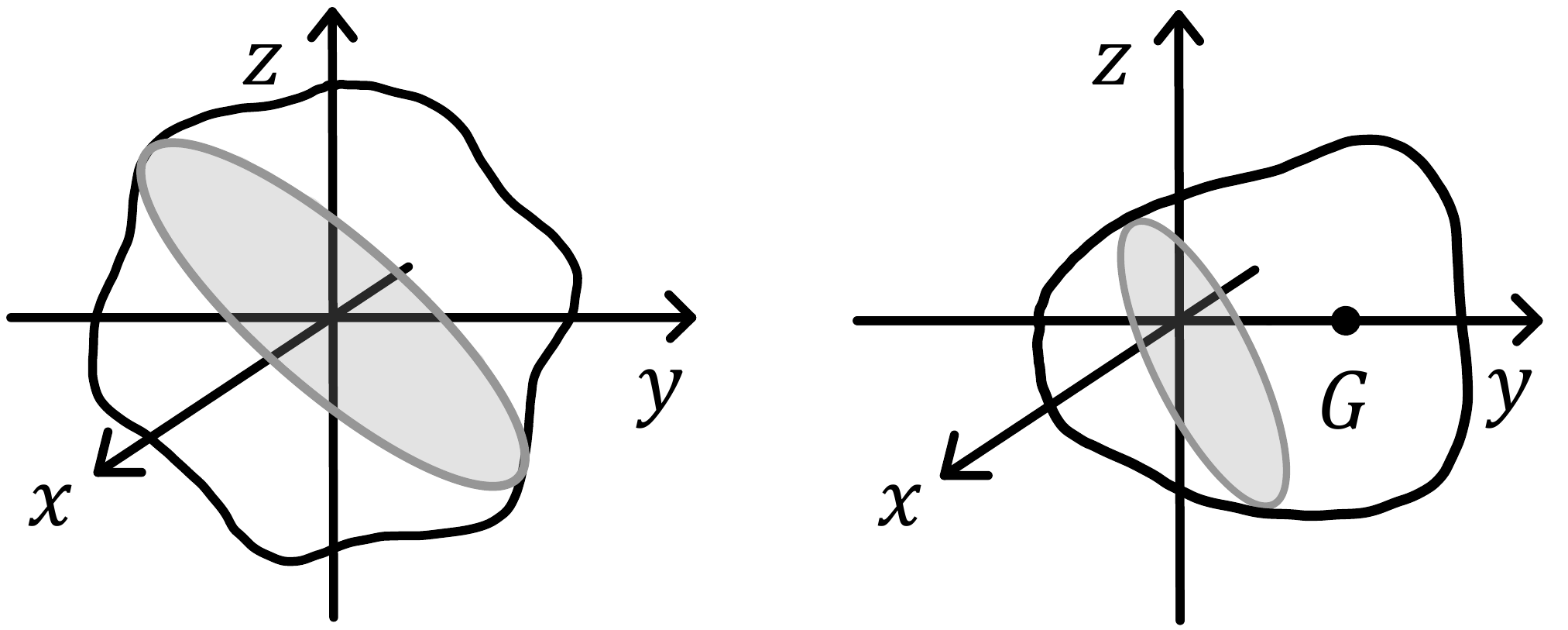}
\caption{\label{fig:1} Comparison of the initial idea and the extended problem in $n=3$, $m=2$. The first picture describes the problem where the centre of mass is the origin. The second picture describes the extended problem where the centre of mass is G.}
\end{figure}
Again, the problem can be extended to maximizing the average $m$-dimensional cross section in the $n$-dimensional shape. While the initial problem is to solve for the $n=2$ and $m=1$, the problem can be extended to higher dimension. For example, when $n=3$ and $m=2$, we are investigating the three-dimensional shape that maximizes its average cross sectional area. The average cross sectional area is defined as the average of all possible two-dimensional sections of the shape through any point inside the shape. The problem can be extended to the dimensions higher than three. Then, the problem has a constraint that the centroid is at $\underline{G}=(G_1,\cdots,G_n)$. It is specified as
\begin{equation}
    \underline{G}=\frac{1}{V}\int\cdots\int\underline{r}\,d_nV=\frac{1}{V}\int\cdots\int\underline{r}\,r^{n-1}\,dr\,d_{n-1}\Omega
\end{equation}
where $V$ is the $n$-dimensional volume of the shape, $d_nV$ is the $n$-dimensional volume element, and $d_{n-1}\Omega$ is $(n-1)$-dimensional solid angle which is introduced at the section 2.1. Or, the coordinates of the centroid is specified as
\begin{align}
\label{eqn:eqlabel}
\begin{split}
G_1=\frac{1}{V}\int\cdots\int \frac{\left(f(\underline{n})\right)^{n+1}}{n+1}\cos{\phi_1}\,d_{n-1}\Omega,
\\
G_2=\frac{1}{V}\int\cdots\int \frac{\left(f(\underline{n})\right)^{n+1}}{n+1}\sin{\phi_1}\cos{\phi_2}\,d_{n-1}\Omega,
\\
\vdots
\\
G_{n-1}=\frac{1}{V}\int\cdots\int \frac{\left(f(\underline{n})\right)^{n+1}}{n+1}\sin{\phi_1}\cdots\sin{\phi_{n-2}}\cos{\phi_{n-1}}\,d_{n-1}\Omega,
\\
G_n=\frac{1}{V}\int\cdots\int \frac{\left(f(\underline{n})\right)^{n+1}}{n+1}\sin{\phi_1}\cdots\sin{\phi_{n-1}}\,d_{n-1}\Omega.
\end{split}
\end{align}
Even though we cannot draw higher dimensional figures explicitly, we can solve the cases through algebraic solutions or numerical analysis.

We derive the Euler-Lagrange equation of the shape assuming that the shape is star shaped, but solving the problem under the assumption is not satisfactory. So we use a different method. We attempt to derive the equation of the stationary shape by examining the shape through making deformations and carrying out the equations from the basic principles of the calculus of variations. This paper explores properties of the equation of the shape, such as the axial symmetry and investigation of the equation in different coordinate systems. We also analyze the characteristics of the shapes with stationary average thickness and study particular cases, such as $(m, n)$ is (1, 2), (1, 3), (2, 3), and (1,5).
\section{The $n$-Dimensional Sphere}
\subsection{Spherical Coordinates in $n$-Dimensions}
Spherical coordinates are advantageous in the presence of the spherical symmetry. Also, concepts such as higher dimensional solid angle are conveniently defined in spherical coordinates. So investigating the spherical coordinate is essential in studying the shapes that exhibits stationary values of average cross section.

The most basic form of spherical coordinates in $n$ dimensions consist of a radial coordinate $r$ and $n-1$ coordinates $\phi_1, \cdots, \phi_{n-1}$ where $\phi_1, \cdots, \phi_{n-2}$ range over $[0, \pi]$ and $\phi_{n-1}$ ranges over $[0,2\pi).$ The $n$-dimensional spherical coordinates $[r, \phi_1, \cdots ,\phi_{n-1}]$ are related to the $n$-dimensional Cartesian coordinates $[x_1,\cdots,x_n]$ by \cite{M}.
\begin{align}
\begin{split}
x_1=r\cos{\phi_1},
\\ x_2=r\sin{\phi_1}\cos{\phi_2},
\\ \vdots
\\ x_{n-1} = r\sin{\phi_1}\cdots \sin{\phi_{n-2}}\cos{\phi_{n-1}},
\\ x_{n} = r\sin{\phi_1}\cdots \sin{\phi_{n-1}}.
\end{split}
\end{align}
These imply that the Cartesian coordinate satisfies
\begin{equation}
    {x_1}^2+\cdots+{x_n}^2=r^2
\end{equation}
and covers all the $\mathbb{R}^n$ space once. Using the spherical coordinate, the $n$-dimensional volume element is given as
\begin{align}
\begin{split}
    d_n V= \biggl\lvert \det \frac{\partial (x_i)}{\partial (r, \phi_j)}\biggr\rvert dr\,d\phi_1 \cdots d\phi_{n-1} 
    \\ = r^{n-1} \sin^{n-2}{\phi_1}\sin^{n-3}{\phi_2} \cdots \sin{\phi_{n-2}}\,dr\,d\phi_1 \cdots d\phi_{n-1}
\end{split}
\end{align}
and similarly $(n-1)$-dimensional differential of solid angle is defined as
\begin{equation}
    d_{n-1} \Omega = \sin^{n-2}{\phi_1}\sin^{n-3}{\phi_2} \cdots \sin{\phi_{n-2}}\,d\phi_1 \cdots d\phi_{n-1}.
\end{equation}
The volume element and differential of solid angle are related by
\begin{equation}
   d_n V= r^{n-1} dr \, d_n \Omega.
\end{equation}
There is a different way to define the spherical coordinates such that the  coordinate is separated by radial coordinate $r$, two-dimensional spherical coordinate $[\theta, \phi]$ and $n-3$ spherical coordinates $[\alpha_1, \alpha_2, \cdots, \alpha_{n-3}]$. $\theta$ ranges over $[0,\frac{\pi}{2}]$, $\phi$ ranges over $[0,2\pi)$, $\alpha_1, \cdots, \alpha_{n-4}$ range over $[0,\pi]$, and $\alpha_{n-3}$ ranges over $[0,2\pi)$. This spherical coordinate is related to the $n$-dimensional Cartesian coordinates $[x_1, \cdots , x_n]$ as
\begin{align}
    \begin{split}
        x_1 = r \cos{\alpha_1} \cos{\theta}, 
        \\
        x_2 = r \sin{\alpha_1}\cos{\alpha_2} \cos{\theta},
        \\
        \vdots
        \\
        x_{n-3}= r \sin{\alpha_1} \cdots \sin{\alpha_{n-4}}\cos{\alpha_{n-3}} \cos{\theta},
        \\
        x_{n-2}= r \sin{\alpha_1} \cdots \sin{\alpha_{n-3}}\cos{\theta},
        \\
        x_{n-1}=r\sin{\theta} \cos{\phi},
        \\
        x_{n}=r\sin{\theta} \sin{\phi}.
    \end{split}
\end{align}
We can simply express $S^n$ as $[r\cos{\theta} y_1, \cdots,\,r\cos{\theta}y_{n-2},\, r\sin{\theta}\cos{\phi},\, r\sin{\theta}\sin{\phi}]$ where $[y_1,\cdots, y_{n-2}]$ is the $(n-3)$-dimensional Cartesian coordinate constructed by spherical coordinates $[\alpha_1, \alpha_2, \cdots, \alpha_{n-3}]$. The volume element in this coordinate is
\begin{align}
\begin{split}
    d_n V= \biggl\lvert \det \frac{\partial (x_i)}{\partial (r, \alpha_j, \theta, \phi)}\biggr\rvert dr\,d\alpha_1 \cdots d\alpha_{n-3}\,d\theta \, d\phi 
    \\ = r^{n-1} \sin^{n-4}{\alpha_1}\sin^{n-5}{\alpha_2} \cdots \sin{\alpha_{n-4}}
    \cos^{n-3}{\theta} \sin{\theta}
    \,dr\,d\alpha_1 \cdots d\alpha_{n-3}\,d\theta \, d\phi
\end{split}
\end{align}
and the $(n-1)$-dimensional differential of solid angle is defined as 
\begin{equation}
    d_{n-1} \Omega = \sin^{n-4}{\alpha_1}\sin^{n-5}{\alpha_2} \cdots \sin{\alpha_{n-4}}
    \cos^{n-3}{\theta} \sin{\theta}\,d\alpha_1 \cdots d\alpha_{n-3}\,d\theta \, d\phi.
\end{equation}
\subsection{Volume and Surface Area of $n$-Dimensional Sphere}
The volume and the surface area of the $n$-dimensional sphere is crucial in exploring this problem. Formula for the volume of the sphere is widely used. For example, the formula is necessary to evaluate the thickness of arbitrary dimension, or integrate the volume or the moment in the hyper cylindrical coordinate. The volume and surface area of $n$-dimensional sphere can be derived through the recurrence relation and the Gaussian integral. \cite{D}-\cite{G}

In the $n$-dimension Euclidean space, the volume of the sphere is proportional to $n$th power of its radius, $R$, and the surface area of the sphere is proportional to $(n-1)$th power of its radius. So $(n-1)$-dimensional surface area can be expressed as ${S}_{n-1}R^n$ and $n$-dimensional volume can be expressed as ${V}_{n}R^n$, where $S_{n-1}$ and $V_n$ are the surface area and volume of the unit $n$-sphere.

First, the volume and the surface area of the $n$-sphere is obtained through two recurrence relations. The $(n+1)$-ball is the union of concentric $n$ spherical shells, and it is expressed as
\begin{equation}
    {S}_{n}r^n=\frac{d}{dR}{V}_{n+1}R^{n+1}=(n+1){V}_{n+1}R^{n}.
\end{equation}
It is equivalent to
\begin{equation}
    {V}_{n+1}=\frac{{S}_{n}}{n+1}.
\end{equation}
Another relation of the volume and the surface area of the sphere can be deduced from the new spherical coordinates introduced in the equation (2.6). Using the solid angle defined in the equation (2.6), we can derive the surface area of the $(n+1)$-sphere as
\begin{align}
    \begin{split}
    S_{n+1}=\int d_{n} \Omega
    \\=\int\cdots\int \sin^{n-2}{\alpha_1}\sin^{n-3}{\alpha_2} \cdots \sin{\alpha_{n-2}}
    \cos^{n-1}{\theta} \sin{\theta}\,d\alpha_1 \cdots d\alpha_{n-1}\,d\theta \, d\phi
    \\ =\frac{2\pi}{n} \int\cdots\int \sin^{n-2}{\alpha_1}\sin^{n-3}{\alpha_2} \cdots \sin{\alpha_{n-2}}\,d\alpha_1 \cdots d\alpha_{n-1}
    \\=\frac{2\pi}{R^n}\int\cdots\int d_nV.
    \end{split}
\end{align}
So, we can derive another recurrence relation
\begin{equation}
    S_{n+1}=2\pi V_n.
\end{equation}
Combining the equation (2.10) and the equation (2.12) gives the recurrence relation of $V_n$ as
\begin{equation}
    V_{n+2}=2\pi\frac{V_n}{n+2}.
\end{equation}
Define ${V}_{0}=1$ and ${S}_{0}=2$. Then the recurrence relation of $V_n$ holds true at $n=0$. The solution of the recurrence relation is
\begin{equation}
    {V}_{2k}=\frac{{\pi}^k}{k!}
\end{equation}
and
\begin{equation}
    {V}_{2k+1}=\frac{2{(2\pi)}^k}{(2k+1)!!}
\end{equation}
where $k$ is a natural number. Or we can use the gamma function for the general case. The gamma function is defined as
\begin{equation}
    \Gamma(p)=\int\limits_{0}^{\infty}e^{-x}x^{p-1}\,dx
\end{equation}
where $p \in \mathbb{C}$ and $Re(p)>0$. The gamma function has following properties:
\begin{itemize}
    \item $\Gamma (1)=1$.
    \item $\Gamma(n+1)=n\Gamma (n)$.
    \item $\Gamma(n+1)=n!$ if $n\in\mathbb{N}$.
    \item $\Gamma(\frac{1}{2})=\sqrt{\pi}$.
\end{itemize}
Using properties of the gamma function, the equation (2.14) and (2.15) can be written as
\begin{equation}
    {V}_{n}=\frac{{\pi}^{\frac{n}{2}}}{\Gamma{(\frac{n}{2}+1)}}.
\end{equation}
Substituting the equation (2.17) to the equation (2.12), we have
\begin{equation}
    {S}_{n}=\frac{2{\pi}^{\frac{n+1}{2}}}{\Gamma(\frac{n+1}{2})}.
\end{equation}
We can derive the volume and the surface of the $n$-dimensional sphere in completely different way. We will use the Gaussian integral. Solving the recurrence relation and using the characteristics of Gaussian integral give the same result for the volume and the surface area of the $n$-sphere. The Gaussian integral in $n$-dimensions is written as
\begin{equation}
    I=\int \limits_{-\infty}^{\infty} \cdots \int \limits_{-\infty}^{\infty} e^{-(
    {{x}_{1}}^2+{{x}_{2}}^2+\,\cdots\,+{{x}_{n}}^2)}\,d{x}_{1} \cdots d{x}_{n}=(\sqrt{\pi})^n={\pi}^{\frac{n}{2}}.
\end{equation}
The Gaussian integral is also written in spherical coordinates such as
\begin{equation}
    I=\int_{{S}^{n-1}} \int \limits_{0}^{\infty}e^{-{r}^{2}}(r^{n-1}\,d_{n-1}\Omega \, dr)=\int_{{S}^{n-1}}d_{n-1}\Omega\int \limits_{0}^{\infty}e^{-{r}^{2}}r^{n-1}\, dr.
\end{equation}
Substituting $r^2=R$,
\begin{equation}
    I=\int_{{S}^{n-1}}d_{n-1}\Omega\int \limits_{0}^{\infty}\frac{e^{-R}R^{\frac{n}{2}-1}}{2}\, dR = \int_{{S}^{n-1}}d_{n-1}\Omega\frac{1}{2}\,\Gamma{(\frac{n}{2})}.
\end{equation}
Therefore, equating the equation (2.19) and the equation (2.21) gives the surface area of the unit $(n-1)$-sphere as \begin{equation}
   {S}_{n-1}= \int_{{S}^{n-1}}d_{n-1}\Omega=\frac{2{\pi}^{\frac{n}{2}}}{\Gamma(\frac{n}{2})}.
\end{equation}
Substituting this result into relationship of the volume and the surface area in equation (2.10), we get
\begin{equation}
    V_n=\frac{S^{n-1}}{n}=\frac{\pi^{\frac{n}{2}}}{\Gamma(\frac{n}{2}+1)}.
\end{equation}
\section{Average $m$-Dimensional Thickness of $n$-Dimensional Shape}
First, assume that the shapes are star shaped. The definition of star shaped region is to have the unique points for each of the radial directions. It is similar to the definition of the function, hence the shape in $n$-dimension is defined as a function. The radial distance of the boundary of the shape is expressed as a function about the direction from the origin. It is expressed as $r=f({\phi}_{1},{\phi}_{2},\cdots{\phi}_{n-1})=f(\underline{n})$, where ${\phi}_{1},{\phi}_{2},\cdots{\phi}_{n-1}$ are the angular coordinates and $\underline{n}$ is the unit vector pointing the boundary of the shape in the radial direction. Expressing the shape in the function gives the average one-dimensional thickness for $n=2$ as
\begin{equation}
    T(1,2)=\frac{1}{\pi}\int f(\theta)\,d\theta.
\end{equation}
Also, the average one-dimensional thickness for $n=3$ is
\begin{equation}
    T(1,3)=\frac{1}{2\pi}\iint f(\theta,\phi)\sin{\theta}\,d\theta\,d\phi.
\end{equation}]
Solving the average one-dimensional thickness is simple. We only have to integrate twice of the radial distance about the solid angle and average it. However, calculating the average thickness of the cross section with higher dimension is more complicated, because the cross sections pass through different axes. So multiple integrals are needed to calculate the average thickness. We would not be able to calculate the average thickness from the conventional methods, by integrating all the possible cases by one by one. We need a general formula to simply calculate the average thickness of cross sections of higher dimensions. The formula can be inferred from the properties of average thickness in arbitrary dimension.

In the general case, the $m$-dimensional thickness of $n$-dimensional shape should be
\begin{itemize}
\item degree $m$ in $f({\phi}_{1},\cdots{\phi}_{n-1})$. It should give the $m$ volume of the cross section.
\item spherically symmetric. Average thickness should not be dependent on coordinate system.
\item must give ${V}_{m}$ when $f({\phi}_{1},\cdots{\phi}_{n-1})\equiv1$. The shape is just a unit $n$-sphere when the radius is equivalent to 1.
\end{itemize}
These three conditions indicate strongly that the formula of $m$-dimensional thickness in $n$-dimensional shape is
\begin{equation}
    T(m,n)=\frac{{V}_{m}}{{S}_{n-1}}\int \cdots \int {f({\phi}_{1},\cdots{\phi}_{n-1})}^{m} \,{d}_{n-1}\Omega.
\end{equation}
$V_m$ and $S_{n-1}$ is computed from the equation (2.26) and the equation (2.27). For example, from the equation (3.3),
\begin{equation}
    T(2,3)=\frac{1}{4}\iint (f(\underline{n}))^2 \, d_2 \Omega .
\end{equation}
This can be justified through the arguments regarding its coordinate system. The area of the cross section through the pole $(0,0,1)$ and $\phi=0$ and $\phi=\pi$ is
\begin{equation}
    T_{(0,0,1), (1,0,0)}=\int \limits_{0}^{\pi}\frac{1}{2}\left(f(\theta,0)\right)^2\,d\theta+    \int \limits_{0}^{\pi}\frac{1}{2}\left(f(\theta,\pi)\right)^2\,d\theta.
\end{equation}
Then, the average area of section through the pole $(0,0,1)$ is
\begin{align}
    \begin{split}
          T_{(0,0,1)}=\frac{1}{2\pi}\int \limits_{0}^{2\pi} \int \limits_{0}^{\pi}\left(f(\theta,\phi)\right)^2\,d\theta \,d\phi=\frac{1}{2\pi}\iint\frac{\left(f(\theta,\phi)\right)^2}{\sin{\theta}}\,d_2\Omega 
          \\ = \frac{1}{2\pi}\iint \frac{\left(f(\underline{n})\right)^2}{\sqrt{1-{\left(\underline{n}\cdot\underline{{e}_{1}}\right)}^2}}d_2\Omega
    \end{split}
\end{align}
where $\underline{e_1}=(0,0,1)$. This is the coordinate-independent formula of the average cross sectional area through the pole $\underline{e_1}$. By the symmetry, the average area of the cross section through the axis $\underline{e}$ is
\begin{equation}
    T_{\underline{e}}=\frac{1}{2\pi}\iint \frac{\left(f(\underline{n})\right)^2}{\sqrt{1-{\left(\underline{n}\cdot\underline{e}\right)}^2}}d_2\Omega.
\end{equation}
We should average all possible $\underline{e}$ for the average cross sectional area. Therefore,
\begin{equation}
    T(2,3)=\frac{1}{4\pi}\iint\left[
    \frac{1}{2\pi}\iint \frac{\left(f(\underline{n})\right)^2}{\sqrt{1-{\left(\underline{n}\cdot\underline{e}\right)}^2}}d_2\Omega\right]\,d_2\Omega'.
\end{equation}
We can change the order of integration as
\begin{equation}
    T(2,3)=\frac{1}{8\pi^2}\iint \left(f(\underline{n})\right)^2\,d_2\Omega \iint \frac{d_2\Omega'}{{\sqrt{1-{\left(\underline{n}\cdot\underline{e}\right)}^2}}}.
\end{equation}
By choosing a certain coordinate for $\underline{e}$, the second integral can be independent of $\underline{n}$. We can choose an orthonormal coordinate for $\underline{e}$, such as
\begin{equation}
    \underline{e}=\left(\underline{n}\cos{\theta}+\underline{n_1}\sin{\theta}\cos{\phi}+\underline{n_2}\sin{\theta}\sin{\phi}
    \right).
\end{equation}
Using this coordinate system,
\begin{equation}
    \iint \frac{d_2\Omega'}{{\sqrt{1-{\left(\underline{n}\cdot\underline{e}\right)}^2}}}=\iint d\theta \, d\phi
    =2\pi^2.
\end{equation}
And the formula for the average cross section is equivalent to the equation (3.4). Also, the average $m$-dimensional thickness in $n$-dimensional shape can be justified through the same arguments. We should start by averaging the $m$-dimensional section passing through orthogonal set of vectors $( \underline{{e}_{1}},\underline{{e}_{2}},\cdots,\underline{{e}_{m-1}})$. Then the average $m$-dimensional thickness $d$ through these vectors will be 
\begin{equation}
    T_{\underline{e_1},\cdots,\underline{e_{m-1}}}=\int\cdots\int\frac{\left(f(\underline{n})\right)^m}
    {[
    1-(\underline{e_1}\cdot\underline{n})^2-(\underline{e_2}\cdot\underline{n})^2-\cdots-(\underline{e_{m-1}}\cdot\underline{n})^2
    ]^{\frac{n-m}{2}}}\,d_n\Omega.
\end{equation}
By the same argument, integrating equation (3.12) by all possible vectors will give an integral independent of $\underline{n}$ by setting appropriate coordinate systems for $( \underline{{e}_{1}},\underline{{e}_{2}},\cdots,\underline{{e}_{m-1}})$. And, again, equation (3.3) can be obtained. 

The average $m$-dimensional thickness in $n$-dimensional shape can be expressed differently. From equation (3.3),
\begin{align}
\begin{split}
     T(m,n)=\frac{{V}_{m}}{{S}_{n-1}}\int \cdots \int {f({\phi}_{1},\cdots{\phi}_{n-1})}^{m} \,{d}_{n-1}\Omega
     \\
     =\frac{{V}_{m}}{{S}_{n-1}}\int \cdots \int mr^{m-1}\,dr
     \,{d}_{n-1}\Omega
     \\
     =\frac{m{V}_{m}}{{S}_{n-1}}\int \cdots \int \frac{r^{n-1}\,dr
     \,{d}_{n-1}\Omega}{r^{n-m}}
     \\
     =\frac{m{V}_{m}}{{S}_{n-1}}\int \cdots \int \frac{d_nV}{r^{n-m}}.
\end{split}
\end{align}
Or, in differential form, the thickness can be written again as
\begin{equation}
    dT(m,n)=\frac{m{V}_{m}}{{S}_{n-1}} \frac{d_nV}{r^{n-m}}.
\end{equation}
The equation (3.13) shows the linearity of the contribution to the average thickness of the each volume elements. This equation is advantageous that it can be applied into the general shapes. Because we expressed the shapes as the function with the assumption that the shapes are star shaped, the equation (3.3) cannot be applied to the shapes which are not star shaped. However, the equation (3.13) can evaluate the average $m$-dimensional thickness of any shapes.
\section{Derivation Through Calculus of Variation}
\subsection{Derivation Through Lagrange Multiplier Method}
To solve the equation for the shape that maximizes the average $m$-dimensional thickness, we solve the Euler-Lagrange equation using the calculus of variation. Because the shapes are considered as a function, $r=f{(\underline{n})}$, only the star shaped regions are solved through the Euler-Lagrange equations. Because of the assumption of the shapes being star shaped, this method does not actually prove what we want to prove. However, this method gives a correct solution, and it is worthwhile to examine this method.

To include the information of the constraints on its volume and centre of mass, Lagrange multiplier method is used. The functional of average thickness of $m$-dimensional cross section in $n$-dimensional shape is written in equation (3.3). And the constraint that its volume should be constant as $V$ is expressed as
\begin{equation}
    V=\int \cdots \int \, d_nV=\int \cdots\int r^{n-1}\,dr\,d_{n-1}\Omega=\int\cdots\int \frac{\left(f(\underline{n})\right)^n}{n}\,d_{n-1}\Omega.
\end{equation}
The constraint on the centre of mass is expressed as the equation (1.4). Using the Lagrange multiplier method, the constraints that the volume is constant as $V$ and centroid is at $\underline{G}$ is introduced in Euler-Lagrange equation as
\begin{align}
\begin{split}
    \frac{\partial}{\partial f}\biggl(
    \frac{V_m}{S_{n-1}}\,f^m-\lambda'\,\frac{f^n}{n}-\mu_1'\,\frac{1}{V}\,\frac{f^{n+1}}{n+1}\cos{\phi_1}-\mu_2'\,\frac{1}{V}\,\frac{f^{n+1}}{n+1}\sin{\phi_1}\cos{\phi_2}
    \\
    -\cdots-\mu_n'\,\frac{1}{V}\,\frac{f^{n+1}}{n+1}\sin{\phi_1}\cdots{\sin{\phi_{n-1}}}
    \biggr)=0.
\end{split}
\end{align}
Differentiation gives 
\begin{align}
\begin{split}
    m\frac{V_m}{S_{n-1}}\,f^{m-1}-\lambda'f^{n-1}-\mu_1'\frac{f^n\cos{\phi_1}}{V}-\mu_2'\frac{f^n\sin{\phi_1}\cos{\phi_2}}{V}
    \\
    -\cdots-\mu_n'\frac{f^n\sin{\phi_1}\cdots\sin{\phi_{n-1}}}{V}=0.
\end{split}
\end{align}
Let $r=f(\underline{n})$ and $\lambda=\frac{S_{n-1}}{mV_m}\lambda'$, $\mu_1=\frac{S_{n-1}}{mV_m}\frac{\mu_1'}{V}$, $\mu_2=\frac{S_{n-1}}{mV_m}\frac{\mu_2'}{V}$, $\cdots$, $\mu_n=\frac{S_{n-1}}{mV_m}\frac{\mu_n'}{V}$. Substituting these gives a simpler form of Euler-Lagrange equation as
\begin{align}
\begin{split}
    r^{m-1}-\lambda r^{n-1}
    -r^n\bigl(\mu_1\cos{\phi_1}-\mu_2\sin{\phi_1}\cos{\phi_2}
    \\
    -\cdots-\mu_n\sin{\phi_1}\cdots\sin{\phi_{n-1}}
    \bigr)=0.
\end{split}
\end{align}
\subsection{Stationary Shape}
\begin{figure}[tbp]
\centering 
\includegraphics[width=.5\textwidth]{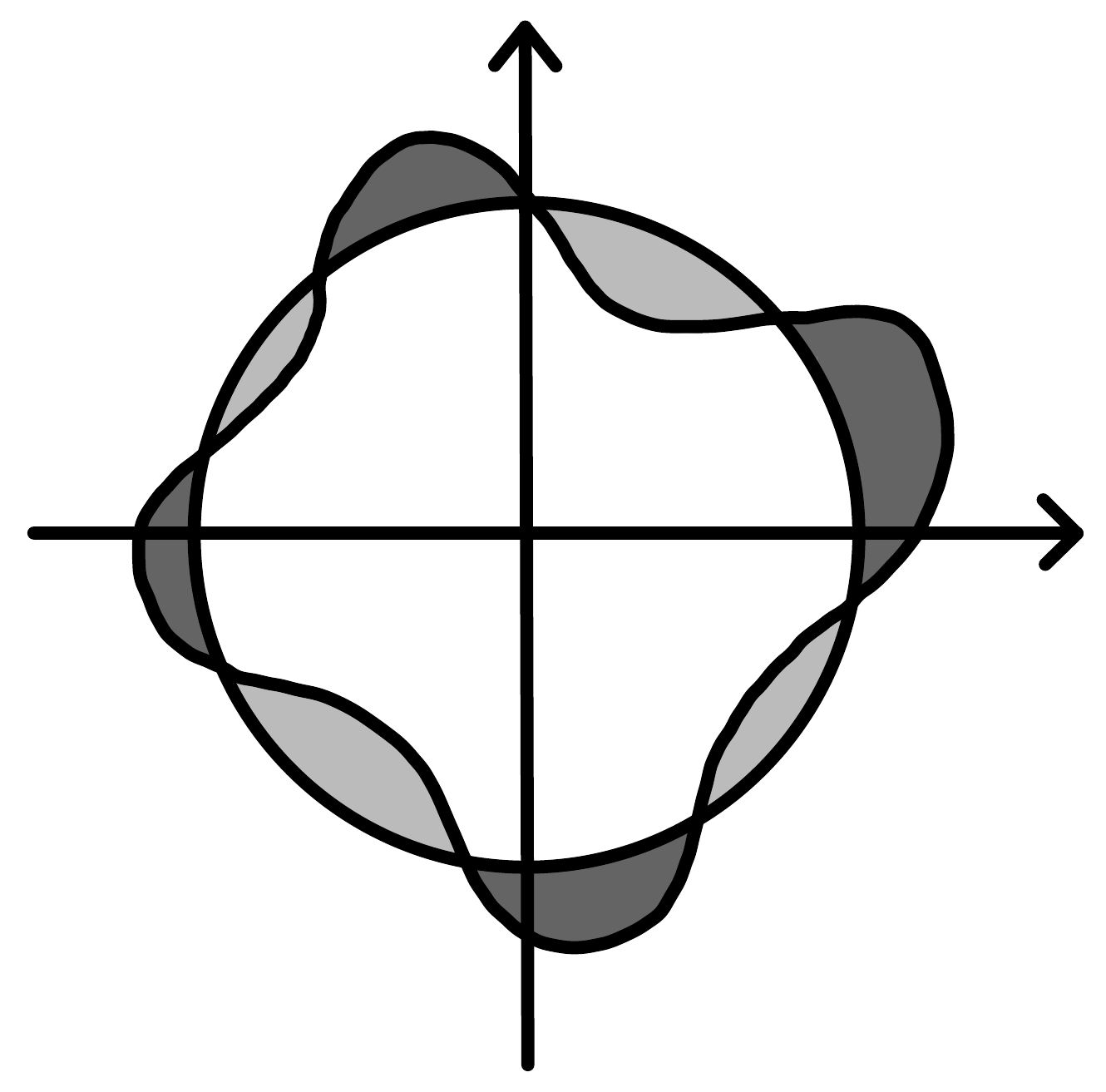}
\caption{\label{fig:2} Any deformation from a sphere decreases the average cross sectional thickness.}
\end{figure}
The solution of the Euler-Lagrange equation does not always give the function that maximizes the average thickness of the $m$-dimensional cross section. Generally, the solution gives the stationary value for the average thickness. However, when the centroid is at the origin, the sphere maximizes the average cross sectional thickness.

According to equation (3.14), the contribution of a volume element to the average thickness is inversely proportional to $r^{n-m}$. In other words, the contribution of a volume element to the average thickness increases as its distance from the origin decreases. Therefore, the average cross section through the centre of mass is maximized when the shape is $n$-dimensional sphere. Any variations on the sphere will decrease the average $m$-dimensional thickness. Figure 2 shows that any deformation from the sphere draws the volume elements out of the sphere, and the average thickness decreases. We can integrate the contribution to the average thickness by each volume elements and explicitly show that the average thickness decreased due to the deformation.
\begin{align}
\begin{split}
      \int \cdots \int_{Outside} \frac{mV_m}{S_{n-1}}\frac{d_nV}{r^{n-m}} \leq \int \cdots \int_{Outside} \frac{mV_m}{S_{n-1}}\frac{d_nV}{R^{n-m}}
      \\
      =\int \cdots \int_{Inside} \frac{mV_m}{S_{n-1}}\frac{d_nV}{R^{n-m}}
      \\
      \leq \int \cdots \int_{Inside} \frac{mV_m}{S_{n-1}}\frac{d_nV}{r^{n-m}}
\end{split} 
\end{align}
where $R$ is the radius of the sphere before the deformation, $Outside$ is the dark gray colored region, and $Inside$ is the light gray colored region.

However, the average cross sectional thickness is not maximized if the centroid is not origin. We produce a proof by considering a dumbbell like shape depicted in Figure 3. The dumbbell like shape is consisted of two spheres and a thin ling neck part that connects two spheres. The dumbbell like shape is actually a slight deformation from the sphere. We move the outermost elements of the sphere far away from the origin to satisfy the constraint on centroid. Because we do not constrain the shape to be the connected shape, we can eliminate the neck part. Therefore, we can simplify this dumbbell like shape as the shape with two $n$-spheres, one with nearly volume $V$ located at the origin, and other one with infinitesimal volume located in the axis that contains the centre of mass, far away from the origin. If we want to include the restriction that the shape is connected, we can make the connection between the two spheres to be infinitesimally thin.

\begin{figure}[tbp]
\centering 
\includegraphics[width=.5\textwidth]{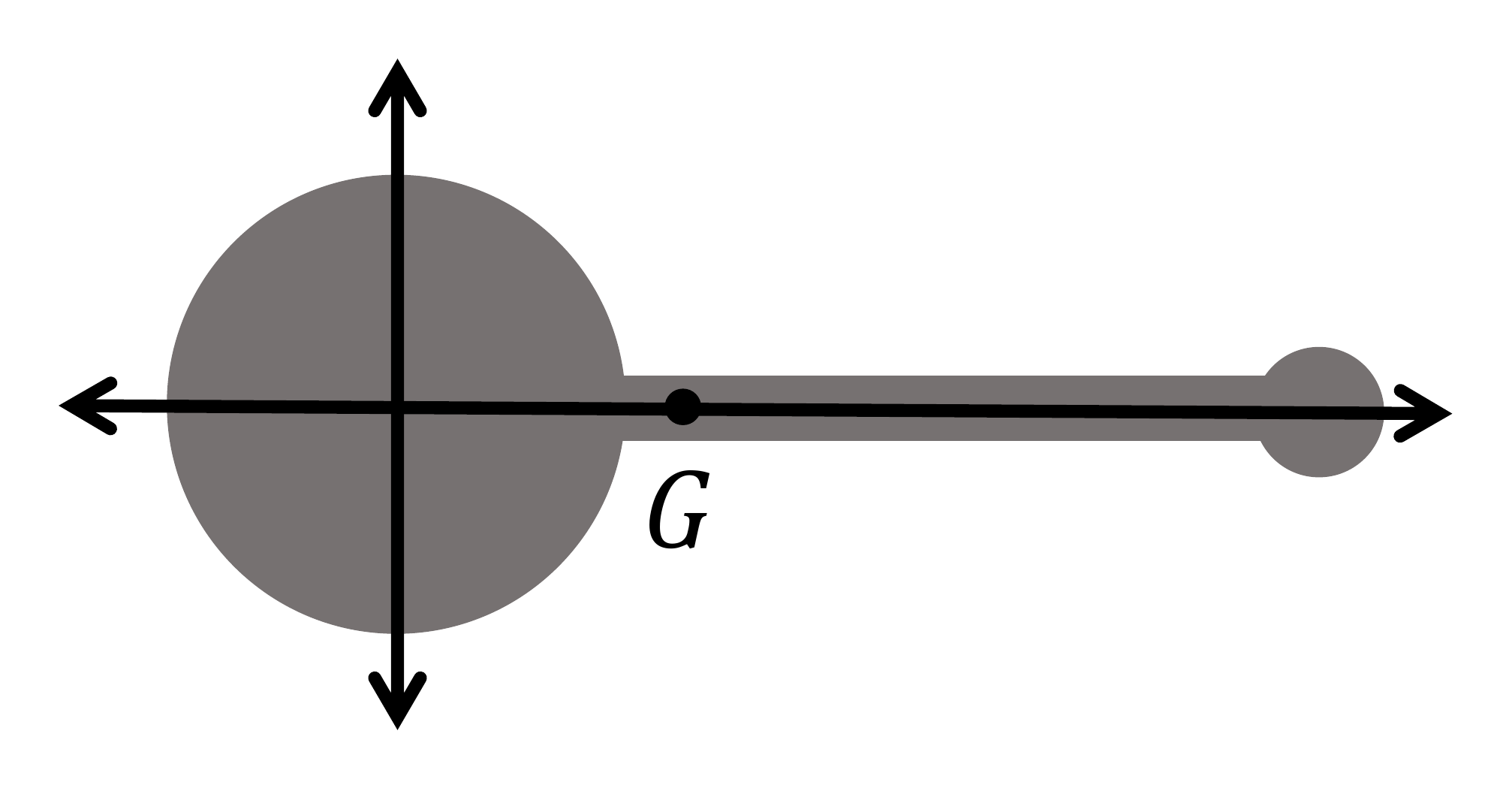}
\caption{\label{fig:3} Dumbbell like shape.}
\end{figure}

Consider the case of the two-dimensional shape with total area $A$ and the centroid fixed at $G$. We can show that this dumbbell like shape can have average diameter greater than $2\sqrt{\frac{A}{\pi}}-\epsilon$, for any $\epsilon > 0$. However, it can never have the average diameter of $2\sqrt{\frac{A}{\pi}}$ unless the centroid is at the origin. This proves that the average diameter of the shape cannot be maximized.

We calculate the average thickness, $T$, of the dumbbell like shape. Let the bigger sphere of the dumbbell have the area $A_1$ and infinitesimal sphere have the area $A_2$, and placed at $x_2$ away from the origin. We assume $A_1 \gg A_2$ and $x_2 \gg 1$. From the constraint on area and centre of mass,
\begin{equation}
    A= A_1+A_2
\end{equation}
and
\begin{equation}
    G=\frac{A_2 x_2}{A}.
\end{equation}
Let the contribution to the average thickness of $A_1$ and $A_2$ be $T_1$ and $T_2$. We can approximate $T_2$ through the equation (3.14) as
\begin{equation}
    T_2=\int\frac{V_1}{S_1}\frac{dA}{r}\approx \frac{V_1}{S_1}\frac{A_2}{x_2}=\frac{A_2}{\pi x_2}.
\end{equation}
Hence, we get
\begin{equation}
    T=T_1+T_2=2\sqrt{\frac{A_1}{\pi}}+\frac{A_2}{\pi x_2}.
\end{equation}
Substituting the equations about the constraints, we have
\begin{equation}
    T=2\sqrt{\frac{A-A_2}{\pi}}+\frac{{A_2}^2}{\pi AG}.
\end{equation}
We can define the ratio of the area of the infinitesimal sphere to the total area as $\gamma = \frac{A_2}{A}$, such that
\begin{align}
    \begin{split}
            T=2\sqrt{\frac{A}{\pi}(1-\gamma)}+\frac{A}{\pi G}\gamma^2.
    \end{split}
\end{align}
We assume $1 \gg \gamma$, and we can approximate the equation (4.11) as
\begin{equation}
    T=2\sqrt{\frac{A}{\pi}}- \sqrt{\frac{A}{\pi}}\gamma+\left(\frac{A}{\pi G}-\frac{1}{4}\sqrt{\frac{A}{\pi}}
\right)\gamma^2+\mathcal{O}(\gamma^3)
\approx 2\sqrt{\frac{A}{\pi}}- \sqrt{\frac{A}{\pi}}\gamma.
\end{equation}
Therefore, if we choose $\gamma < \sqrt{\frac{\pi}{A}}\epsilon$ and $A_2<A\sqrt{\frac{\pi}{A}}\epsilon$, the average diameter is greater than $2\sqrt{\frac{A}{\pi}}-\epsilon$, for any $\epsilon>0$. As the infinitesimal sphere is located further from the origin and has reduced area, the average diameter gets closer to $2\sqrt{\frac{A}{\pi}}$, but not the same, unless the centre of mass is at its origin. Therefore, the maximum cross sectional thickness does not exist.

The shapes obtained through Euler-Lagrange equations are not dumbbell like shape. The shapes do not maximize the average thickness. Instead, the shapes behave like a local maximum or a stationary point. The average thickness decreases as the consequence of slight deformation. Therefore, the Euler-Lagrange equation gives the stationary shape, the shape with stationary average value of thickness.

\section{Derivation Through Deformation Method}
The equation of the shape can be derived through the deformation method. This method applies the basic concept of calculus of variation to this problem; we give a slight variation to any stationary shapes. While the Euler-Lagrange equation is only able to consider the star shaped regions, this method can examine the properties of the non star shaped regions. Also, this problem does not contain any derivatives, so it is not necessary to solve the Euler-Lagrange equation.

First, we assume that the shape is stationary, and investigate if slight deformations do not change its average thickness of $m$-dimensional cross sections. In other words, we examine the property of the stationary shape, about which conditions it should satisfy to be stationary. For stationary $n$-dimensional shape, select any $n+2$ points on the boundary of the shape and deform by 
\begin{equation}
    \delta V_i={r_i}^{n-1}\, \delta r \, \delta_{n-1} \Omega
\end{equation}
at $(r_i,\phi_{1i},\phi_{2i},\cdots,\phi_{(n-1)i})$, for $i=1,\cdots,n+2$. The deformed volume can be either positive or negative. Still, the deformed shape is constrained to have consistent volume, as $V$. With the assumption that the shape is stationary, the shape should satisfy $n+2$ constraints.
\begin{enumerate}
    \item The shape is stationary, so the average thickness should be same after the deformation.
    \begin{equation}
         T(m,n)
     =\frac{m{V}_{m}}{{S}_{n-1}}\int \cdots \int r^{m-1}\,dr
     \,{d}_{n-1}\Omega
    \end{equation}
    should be consistent. Or in discrete version,
    \begin{equation}
        \sum_{i=1}^{n+2} {r_{i}}^{m-1} \,\delta r \, \delta_{n-1} \Omega=0.
    \end{equation}
    \item The volume is constrained to be constant. So the sum of the deformed volume is set to be 0. Therefore,
    \begin{equation}
        \sum_{i=1}^{n+2} \delta V_i = \sum_{i=1}^{n+2}{r_i}^{n-1}\, \delta r \, \delta_{n-1} \Omega=0.
    \end{equation}
    \item There are the rest $n$ constraints on the centre of mass. The centre of mass is fixed at $\underline{G}$, or the moment $\underline{M}=\int \underline{r} \, d_nV$ is fixed. In spherical coordinates, the constraint can be written again as
\begin{align}
\begin{split}
    \sum_{i=1}^{n+2} x_1 \, \delta V_i = \sum_{i=1}^{n+2} {r_i}^n \cos{\phi_1} \, \delta r \, \delta_{n-1} \Omega = 0 
    \\
    \sum_{i=1}^{n+2} x_2 \, \delta V_i = \sum_{i=1}^{n+2} {r_i}^n \sin{\phi_1}\cos{\phi_2} \, \delta r \, \delta_{n-1} \Omega = 0
    \\
    \vdots
    \\
    \sum_{i=1}^{n+2} x_n \, \delta V_i = \sum_{i=1}^{n+2} {r_i}^n \sin{\phi_1}\cdots\sin{\phi_{n-1}} \, \delta r \, \delta_{n-1} \Omega = 0.
\end{split}
\end{align}
\end{enumerate}
\begin{figure}[tbp]
\centering 
\includegraphics[width=.45\textwidth]{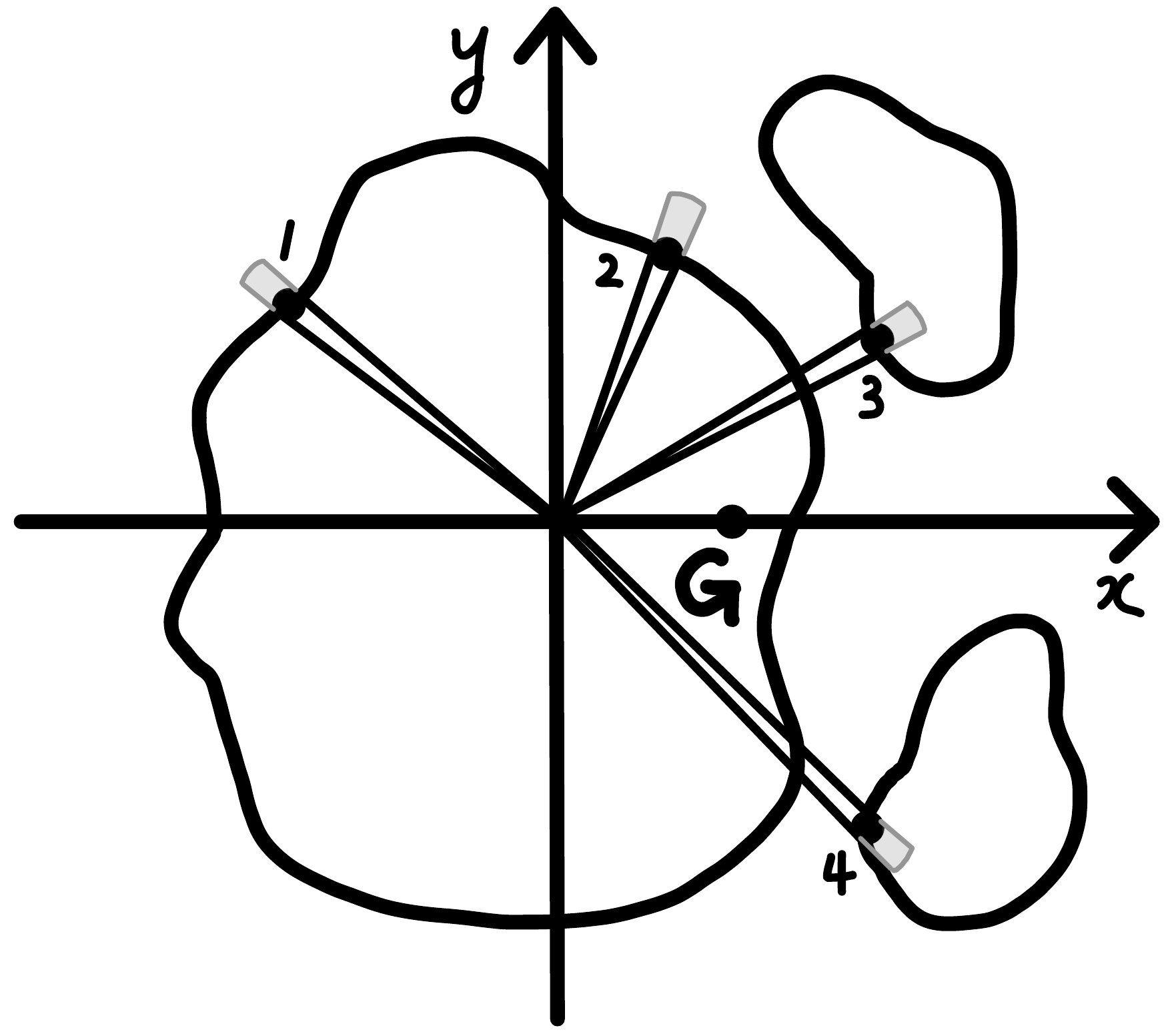}
\caption{\label{fig:4} Deformation in $n=2$ and $m=1$.}
\end{figure}
$\delta r\, \delta_{n-1} \Omega$ satisfies these homogeneous equations. So, the determinant of the matrix $R$ representing the homogeneous equation is 0. The matrix is
\begin{equation}
    R=
    \begin{pmatrix}
  {r_1}^{m-1} & {r_1}^{n-1} & {r_1}^{n}\cos{\phi_1} & {r_1}^{n}\sin{\phi_1}\cos{\phi_2} & \cdots & {r_1}^{n}\sin{\phi_1}\cdots\sin{\phi_{n-1}} \\
  {r_2}^{m-1} & {r_2}^{n-1} & {r_2}^{n}\cos{\phi_1} & {r_2}^{n}\sin{\phi_1}\cos{\phi_2} & \cdots & {r_2}^{n}\sin{\phi_1}\cdots\sin{\phi_{n-1}} \\
  \vdots  & \vdots  & \vdots  & \vdots  & \ddots & \vdots  \\
  {r_n}^{m-1} & {r_n}^{n-1} & {r_n}^{n}\cos{\phi_1} & {r_n}^{n}\sin{\phi_1}\cos{\phi_2} & \cdots & {r_n}^{n}\sin{\phi_1}\cdots\sin{\phi_{n-1}} \\ 
 \end{pmatrix}.
\end{equation}
The rank of the matrix is less than $n+2$, or there is a nonzero vector $(1, \, \lambda, \, \mu_1, \, \mu_2, \,$ $ \cdots, \, \mu_n)$ such that 
\begin{align}
\begin{split}
    r_i^{m-1}-\lambda r_i^{n-1}
    -r_i^n\mu_1\cos{\phi_1}-r_i^n\mu_2\sin{\phi_1}\cos{\phi_2}-\cdots
    \\
    -r_i^n\mu_n\sin{\phi_1}\cdots\sin{\phi_{n-1}}
    =0.
\end{split}
\end{align}
The equation does not depend on which $n+2$ points we select. Therefore, the equation should be true for all points on the boundary shape. Equation (4.6) is obtained again by using different method.
\section{Euler-Lagrange Equation of the Shape}
\subsection{Axial Symmetry}
Even though equation (4.6) contains the full information of the stationary shape, the equation can not be solved analytically in the current form. The equation can be simplified by the argument on axial symmetry of the shape. The stationary shapes are symmetric about the axis that contains the origin and the centre of mass. This can be deduced from the equation (3.13). The contribution to the average $m$-dimensional thickness by a volume element is inversely proportional to $r^{n-m}$. Any variation to the non symmetric shapes to construct the symmetric shape will bring the volume elements closer to the origin and thereby increase the average thickness. So the stationary shapes should be symmetric about its axis.
\begin{figure}[tbp]
\centering 
\includegraphics[width=.6\textwidth]{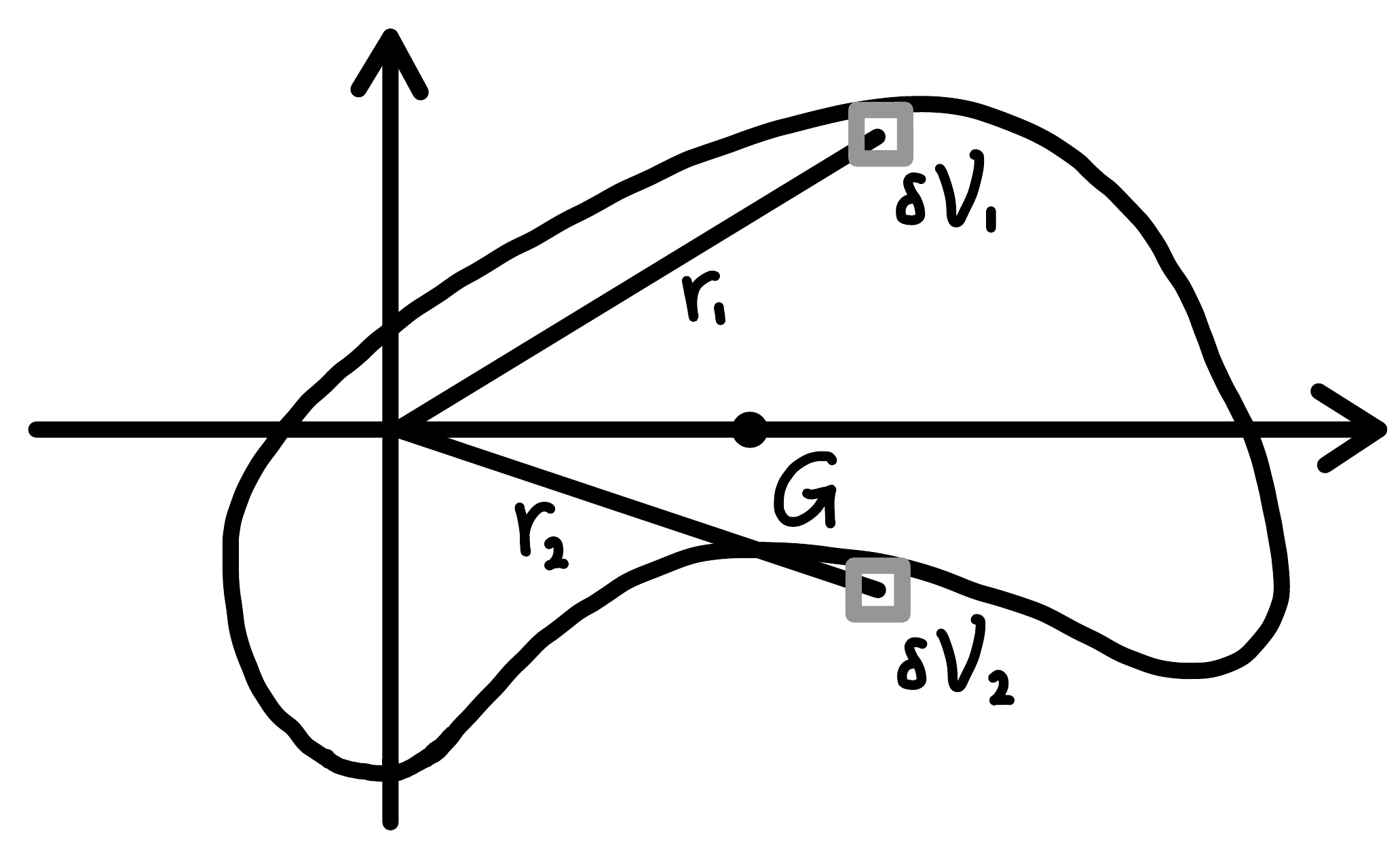}
\caption{\label{fig:5} Moving $\delta V_1$ to $\delta V_2$ increases the average thickness of the shape.}
\end{figure}
With the axial symmetry of the shape, equation (4.6) can be simplified. We are going to take the polar axis as the axis that contains the centre of mass. First, define $\mu={\mu_1}^2+{\mu_2}^2+\cdots+{\mu_n}^2$. And define two unit vectors as
\begin{equation}
    \underline{e}=\frac{1}{\sqrt{{\mu_1}^2+{\mu_2}^2+\cdots+{\mu_n}^2}}\left({\mu_1},{\mu_2},\cdots,{\mu_n}\right)=\frac{1}{\mu}\left({\mu_1},{\mu_2},\cdots,{\mu_n}\right)
\end{equation}
and
\begin{equation}
    \underline{v}=\left(\cos{\phi_1},\,\sin{\phi_1}\cos{\phi_2},\,\cdots,\,\sin{\phi_1}\cdots\sin{\phi_{n-1}}\right).
\end{equation}
$\underline{e}$ is parallel to the axis containing the centre of mass and $\underline{v}$ is the unit vector pointing the boundary of the shape in the radial direction. The equation (4.6) can be written again as
\begin{align}
\begin{split}
    r^{m-1}-\lambda r^{n-1}
    -r^n\bigl(\mu_1\cos{\phi_1}-\mu_2\sin{\phi_1}\cos{\phi_2}
    \\
    -\cdots-\mu_n\sin{\phi_1}\cdots\sin{\phi_{n-1}}
    \bigr)
\\
    =r^{m-1}-\lambda r^{n-1}-\mu r^n\left(\underline{e}\cdot\underline{v} \right)
    \\
    =r^{m-1}-\lambda r^{n-1}-\mu r^n\cos{\theta}=0.
\end{split}
\end{align}
$\theta$ is the polar angle, when the polar axis is set as the axis containing the centre of mass. Dividing equation (6.3) by $r^{m-1}$ gives
\begin{equation}
    1-\lambda r^{n-m} -\mu r^{n-m+1} \cos{\theta} = 0.
\end{equation}
Equation (6.4) shows the axial symmetric nature of the stationary shapes.
\subsection{Cylindrical Coordinates}
Due to the axial symmetric property of the stationary shapes, hyper cylindrical coordinate is advantageous to describe the shapes. The coordinate is separated by the $n-1$ spherical coordinates $[R,\phi_1,\cdots,\,\phi_{n-2}]$ and an axial coordinate $z$, where $R$ ranges over $[0,\infty)$, $\phi_1,\cdots,\phi_{n-3}$ range over $[0,\pi]$, $\phi_{n-2}$ ranges over $[0,2\pi)$, and $z$ ranges over $(-\infty,\infty)$. This hyper cylindrical coordinate $[R,\phi_1,\cdots,\,\phi_{n-2}]$ is related to the $n$-dimensional Cartesian coordinates $[x_1,\cdots,x_n]$ as
\begin{align}
    \begin{split}
        x_1=R\cos{\phi_1},
        \\
        \vdots
        \\
        x_{n-2}=R\sin{\phi_1}\cdots \sin{\phi_{n-3}} \cos{\phi_{n-2}},
        \\
        x_{n-1}=R\sin{\phi_1}\cdots \sin{\phi_{n-2}},
        \\
        x_n=z.
    \end{split}
\end{align}
It implies that
\begin{equation}
    R^2+z^2=r^2
\end{equation}
where $r$ is the radial distance defined in $n$-dimensional spherical coordinate. Therefore, the equation (6.4) can be written in cylindrical coordinate as
\begin{align}
\begin{split}
    1-\lambda r^{n-m}-\mu r^{n-m+1}\cos{\theta}=1-\lambda r^{n-m}-\mu r^{n-m}z
    \\
    =1-r^{n-m}(\lambda+\mu z)=0.
\end{split}
\end{align}
Substituting equation (6.6) in equation (6.7) gives
\begin{equation}
    R^2=\left(\frac{1}{\lambda+\mu z}\right)^{\frac{2}{n-m}}-z^2.
\end{equation}
\section{Stationary Shapes}
\subsection{General Discussions About The Stationary Shapes}
\begin{figure}[tbp]
\centering 
\includegraphics[width=.9\textwidth]{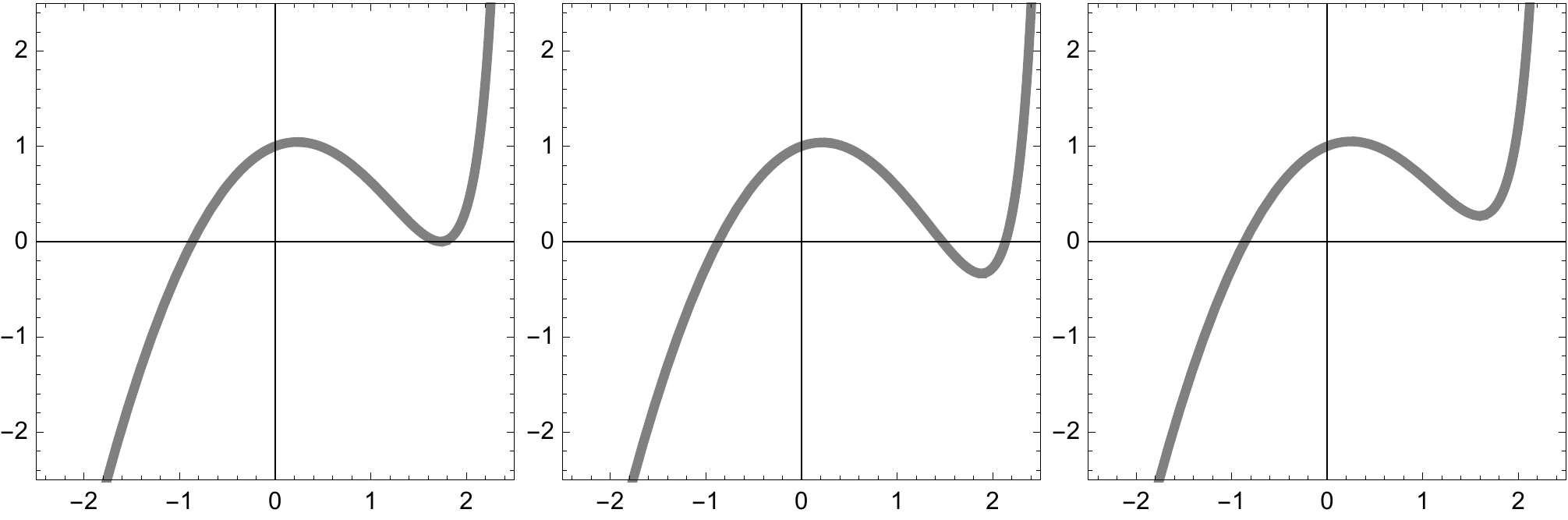}
\caption{\label{fig:6} Three possible cases of $R^2$.}
\end{figure}
\begin{figure}[tbp]
\centering 
\includegraphics[width=.9\textwidth]{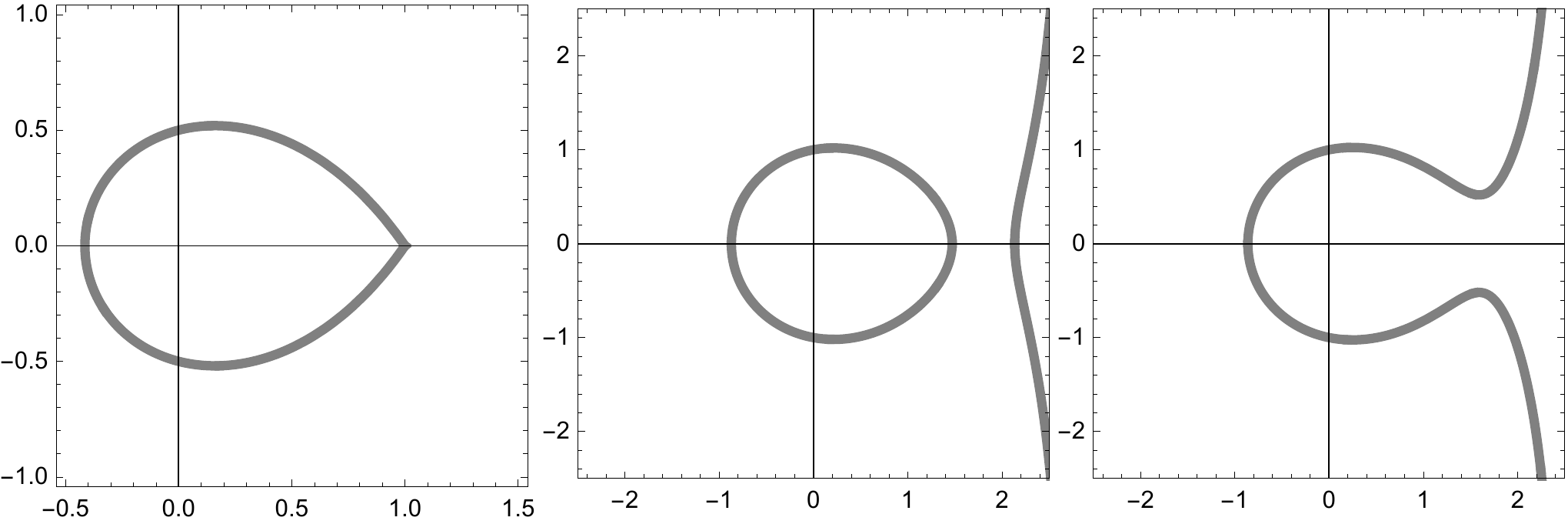}
\caption{\label{fig:7} Three possible cases of $R(z)=R$.}
\end{figure}
Both the spherical equation of the shape, (6.4), and the cylindrical equation of the shape, (6.8), are only dependent on $n-m$. The equation of the stationary shape only depends on the difference of the dimension between the shape and the cross section. Therefore, the spherical equation of the shape of $n-m>4$ cannot be solved algebraically. Yet, the general features of the stationary shape can be examined through cylindrical equation. Generally, equation (6.8) has three possible cases as plotted in Figure 6.
\begin{enumerate}
    \item $R^2$ curve touches the $z$ axis at $z>0$.
    \item $R^2$ curve passes through the $z$ axis twice at $z>0$.
    \item $R^2$ curve only passes through $z$ axis at $z<0$.
\end{enumerate}
Because the cylindrical equation of the stationary shape is actually the square root of equation (6.8), there are also three possible cases for the stationary regions depicted in the hyper cylindrical coordinates. Figure 7 shows three possible cases of $R=R(z)$. First two cases are shapes, but the third case is not actually defined as the stationary shape that we are looking for. It is an opened region without definite volume. The stationary shape corresponding to the second case is an egg shape, which looks like a slight deformation from a sphere. Between these two cases, there is a sharp critical shape which corresponds to the first case in Figure 7. For the critical shape, $R^2$ has a double root at the intersection.

The type of the shape is determined through relationship between the constants $\mu$, $\lambda$, and the difference of the dimension of the shape and the cross section. The relationship between $\mu$ and $\lambda$ in the critical shape is deduced from $R^2$ having double root. From equation (6.8)
\begin{equation}
    1-\lambda z^{n-m} - \mu z^{n-m+1} = 0
\end{equation}
and
\begin{equation}
    -(n-m)\lambda z^{n-m-1}-\mu(n-m+1)z^{n-m}=0.
\end{equation}
Substituting equation (7.2) in equation (7.3) gives
\begin{equation}
    \mu = - (n-m)\,\frac{\lambda^{\frac{n-m+1}{n-m}}}{{(n-m+1)^{\frac{n-m+1}{n-m}}}}.
\end{equation}
In the general case,
\begin{equation}
\mu = - (n-m)\,\frac{\lambda^{\frac{n-m+1}{n-m}}}{{(n-m+1)^{\frac{n-m+1}{n-m}}}}\,e
\end{equation}
where $e \geq 0$. The type of shape is determined by $e$, which is an analogous to eccentricity of the conic section.
\begin{itemize}
    \item If $e=0$, the stationary shape is the sphere.
    \item If $0<e<1$, the stationary shape is the egg shape.
    \item If $e=1$, the stationary shape is the critical shape.
    \item If $1<e$, the shape is opened, and the volume of the shape is not defined.
\end{itemize}
The volume, the moment, and the average thickness of the shape is obtained through integration, and written in terms of the constants, $\mu$ and $\lambda$. Furthermore, the relationship between the volume, the moment, and average $m$-dimensional thickness is directly obtained from the equation (6.3). Multiplying $r$ and integrating the equation about $d_{n-1} \Omega$ gives
\begin{equation}
    \int r^{m} \,d_{n-1}\Omega-\int \lambda r^{n}\,d_{n-1}\Omega-\int\mu r^{n+1}\cos{\theta}\,d_{n-1}\Omega=0
\end{equation}
which is
\begin{align}
    \begin{split}
    \frac{S_{n-1}}{V_m}\left(\frac{V_m}{S_{n-1}}\int \cdots \int r^m \,d_{n-1}\Omega \right)
    +\lambda n \left( \int \cdots \int \frac{r^n}{n} \,d_{n-1}\Omega \right)
    \\
    +\mu (n+1) \left(\int \cdots \int \frac{r^{n+1}}{n+1} \,d_{n-1}\Omega \right)
    \\
    = \frac{S_{n-1}}{V_m}\, T(m,n)+\lambda n V + \mu (n+1) M =0
    \end{split}
\end{align}
where $T(m,n)$ is average $m$-dimensional thickness, $V$ is volume, and $M$ is the moment of the shape. For the remaining parts of section 7, various cases will be discussed, and the volume, the moment, and the average cross sectional thickness will be calculated in terms of $\lambda$ and $\mu$.

\subsection{n = 2, m = 1}
When $n=2$ and $m=1$, the equation of the stationary shape is
\begin{equation}
    1-\lambda r - \mu r^2 \cos{\theta} = 0.
\end{equation}
Solving the quadratic equation, we obtain
\begin{equation}
    r=\frac{-\lambda + \sqrt{\lambda^2+4\mu\cos{\theta}}}{2\mu \cos{\theta}}=
    \frac{2}{\lambda}\,\frac{1}{1+\sqrt{1+\frac{4\mu}{\lambda^2}\cos{\theta}}}.
\end{equation}
From equation (7.4),
\begin{equation}
    \mu=-\frac{\lambda^2}{4}e.
\end{equation}
Relationship between $\mu$ and $\lambda$ can be directly deduced from the equation (7.8). $\lvert\frac{4\mu}{\lambda^2}\rvert=1$ makes the value inside the square root 0, and the sign is set such that the centre of mass is at $z>0$.
\begin{figure}[tbp]
\centering 
\includegraphics[width=.9\textwidth]{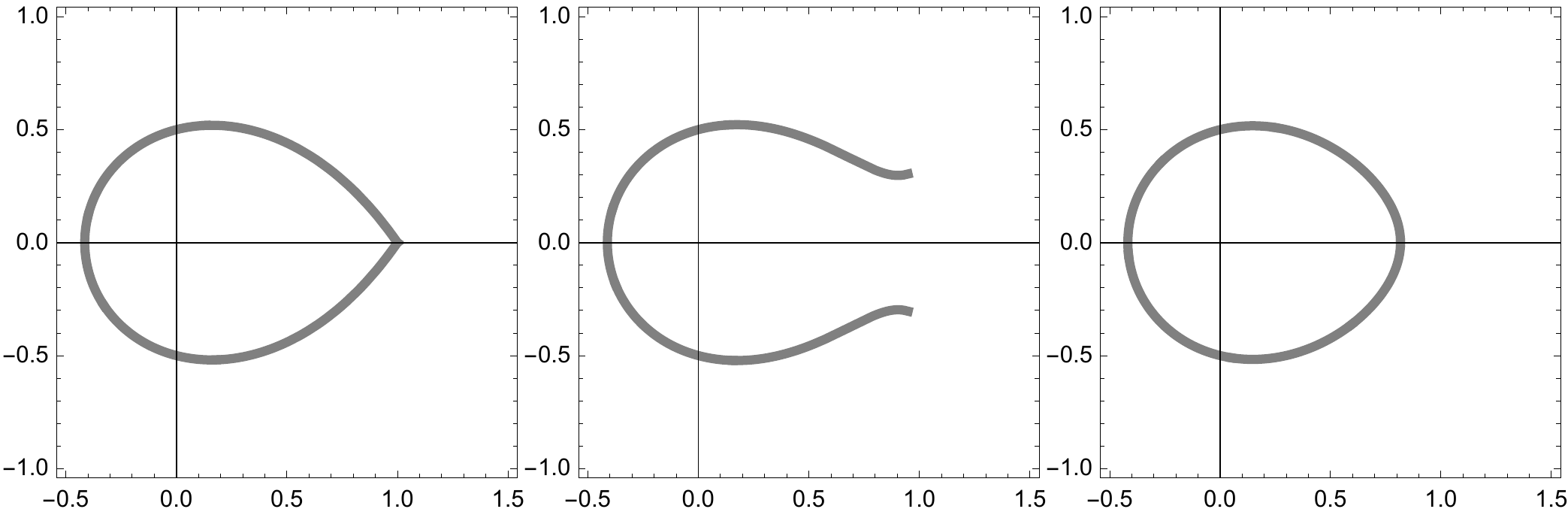}
\caption{\label{fig:8} Stationary regions of $n=2$ and $m=1$. $e=1$, $e=1.05$, and $e=0.95$ from left to right.}
\end{figure}
For general cases, integration does not give the relationship between $\mu$, $\lambda$ and area, moment with the elementary functions. However, the relation can be solved in critical shape. For the critical case, the average thickness $T(1,2)$, area of the shape $A$, and moment $M$ is 
\begin{equation}
    T(1,2)=\frac{2}{\pi\lambda}\int \limits_{0}^{2\pi} \frac{1}{1+\sqrt{1+\cos{\theta}}}\,d\theta=\frac{4\ln{(3+2\sqrt2)}}{\pi\lambda}
\end{equation}
\begin{equation}
    A=\frac{2}{\lambda^2}\int \limits_{0}^{2\pi} \left(\frac{1}{1+\sqrt{1+\cos{\theta}}}\right)^2d\theta=
    \frac{8\sqrt{2}-8\ln{(1+\sqrt{2})}}{\lambda^2}
\end{equation}
\begin{equation}
    M=\frac{4}{\lambda^3}\int \limits_{0}^{2\pi} \left(\frac{1}{1+\sqrt{1+\cos{\theta}}}\right)^3\cos{\theta}\,d\theta
    =\frac{24\ln{(3+2\sqrt{2})}-32\sqrt{2}}{\lambda^3}.
\end{equation}
Furthermore, there is a cusp in the critical shape. The angle of the cusp can be obtained through changing the coordinate system into the Cartesian coordinate. In Cartesian coordinate, the equation (7.7) is written as
\begin{equation}
    \left(\frac{2}{\lambda}\right)^4=\left(\frac{4}{\lambda}+x\right)^2(x^2+y^2).
\end{equation}
From equation (7.13), critical shape in two-dimension is an inverse parabola. Also, the angle of the cusp is $2\arctan{\sqrt{2}}$.
\subsection{n = 3, m = 1}
When $n=3$ and $m=1$, the equation of the stationary shape is 
\begin{equation}
    1-\lambda r^2-\mu r^3 \cos{\theta}=0.
\end{equation}
Hence, equation (7.4) gives
\begin{equation}
    \mu=-\frac{2\sqrt{3}}{9}\lambda^{\frac{3}{2}}e.
\end{equation}
The equation of the shape in spherical coordinate can be obtained by using Cardano method to solve cubic equation. Substitute
\begin{equation}
    \rho = \frac{\sqrt{3}}{\sqrt{\lambda}r}
\end{equation}
in equation (7,14), we get
\begin{equation}
    \rho^3-3\rho+2e\cos{\theta}=0.
\end{equation}
Using Cardano method, let $\rho=\omega+\frac{1}{\omega}$, the equation (7.17) is developed to be
\begin{equation}
    \omega^6+2e\cos{\theta}\omega^3+1=0.
\end{equation}
Quadratic formula gives
\begin{equation}
    \omega^3=-e\cos{\theta}\pm \sqrt{e^2\cos^2{\theta}-1}.
\end{equation}
Define
\begin{equation}
    \cos{\theta_e}=e\cos{\theta}.
\end{equation}
Substituting $\cos{\theta_e}$ gives
\begin{equation}
    \omega=e^{i\left(\frac{\pi-\theta_e}{3}\right)}\,\,or\,\,e^{i\left(\frac{3\pi-\theta_e}{3}\right)} \,\,or\,\, e^{i\left(\frac{5\pi-\theta_e}{3}\right)}.
\end{equation}
From the definition of $\omega$, we obtain the equation of the shape as
\begin{align}
    \begin{split}
        r=\frac{\sqrt{3}}{2\sqrt{\lambda}}\sec{\left( \frac{\pi-\arccos{(e\cos{\theta})}}{3}\right)}
        \\
        or \,\, \frac{\sqrt{3}}{2\sqrt{\lambda}}\sec{\left(\frac{3\pi-\arccos{(e\cos{\theta})}}{3}\right)}
        \\
        or \,\, \frac{\sqrt{3}}{2\sqrt{\lambda}}\sec{\left(\frac{5\pi-\arccos{(e\cos{\theta})}}{3}\right)}.
    \end{split}
\end{align}
\begin{figure}[tbp]
\centering 
\includegraphics[width=.8\textwidth]{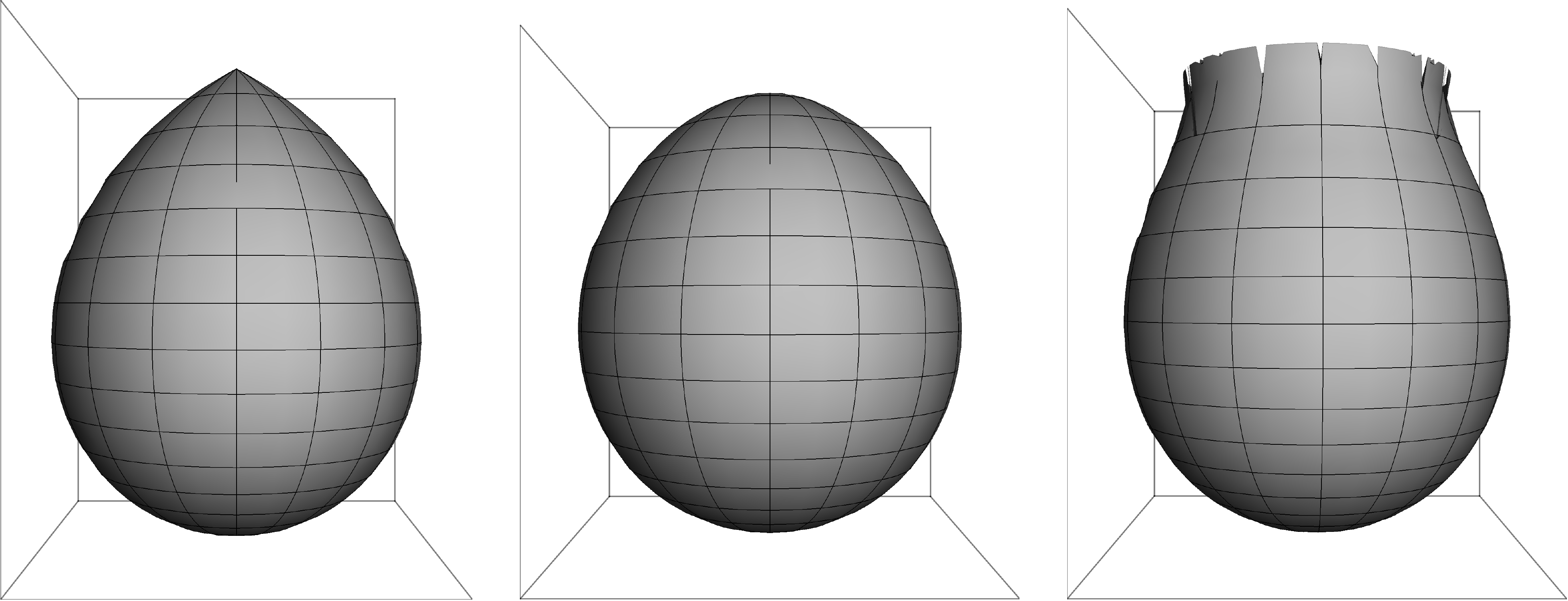}
\caption{\label{fig:9} Stationary regions of $n=3$ and $m=1$. $e=1$, $e=0.95$, and $e=1.1$ from left to right.}
\end{figure}
There are three solutions for the radial distance. This can be explained by the second subfigure of Figure 7. For $\cos{\theta}>0$, or $\theta \in \left(-\frac{\pi}{2},\,\frac{\pi}{2}\right)$, the shape has two positive radial distances and one negative radial distance. And for $\cos{\theta}<0$, or $\theta \in \left(\frac{\pi}{2},\,\frac{3\pi}{2}\right)$, the shape has one positive radial distance and two negative radial distances. So the solution needed to examine the stationary shape is the first solution.

Volume and moment of the shape is integrated in elementary functions only in critical shape. It is easier to integrate in cylindrical coordinate.
\begin{equation}
    V=\pi \int \limits_{z_{-}}^{z_{+}} R^2\,dz=3\sqrt{3}\pi \lambda^{-\frac{3}{2}}\left(\ln{2}-\frac{3}{8}\right)
\end{equation}
and
\begin{equation}
   M= \pi \int \limits_{z_{-}}^{z_{+}} z R^2\,dz=\frac{27}{32}\pi \lambda^{-2} \left( 16\ln2-\frac{21}{2} \right)
\end{equation}
where $z_{+}=\frac{\sqrt{3}}{\sqrt{\lambda}}$ and $z_{-}=-\frac{\sqrt{3}}{2\sqrt{\lambda}}$. Average thickness, which is calculated through the equation (7.6), is
\begin{equation}
    T(1,3)=\frac{3\sqrt{3}}{4\sqrt{\lambda}}(3-2\ln{2}).
\end{equation}
\subsection{n = 3, m = 2}
The equation of the stationary shapes for $n=3$, $m=2$ is equivalent to the equation of the stationary shapes for $n=2$, $m=1$. Let $a=\frac{2}{\lambda}$, the equation (7.8) can be written again as
\begin{equation}
    r=\frac{a}{1+\sqrt{1-e\cos{\theta}}}.
\end{equation}
\begin{figure}[tbp]
\centering 
\includegraphics[width=.8\textwidth]{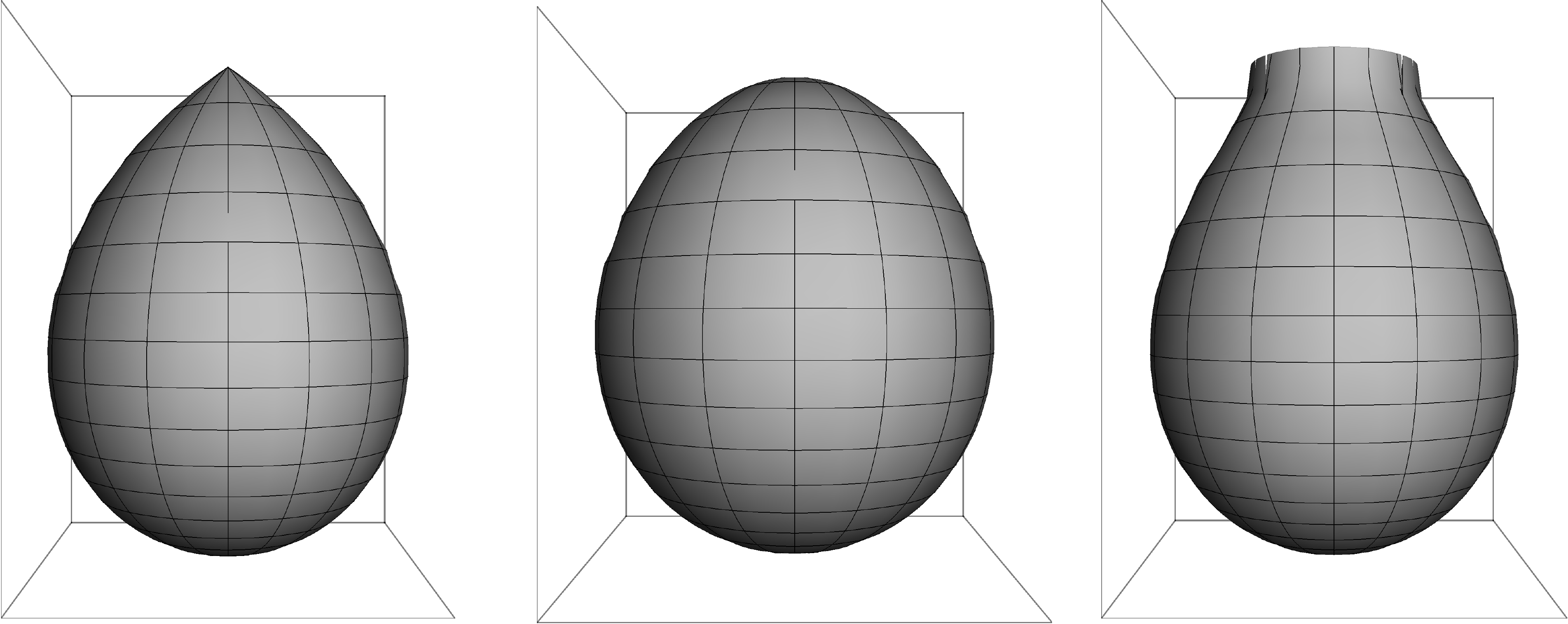}
\caption{\label{fig:10} Stationary regions of $n=3$ and $m=2$. $e=1$, $e=0.95$, and $e=1.03$ from left to right.}
\end{figure}
$n=3$, $m=2$ seems more complicated than $n=2$, $m=1$, but this is actually the simplest case among all cases of maximizing cross section problem. In $n=3$, $m=2$, average cross sectional thickness, volume, and moment can be integrated into elementary functions, because of the $\sin{\theta}$ term in two-dimensional solid angle. We can perform the integration by substitution by substituting $\cos{\theta}$ into another variable. The volume of the shape is
\begin{align}
    \begin{split}
    V = \frac{2\pi}{3} \int \limits_{0}^{\pi} \left( \frac{a}{1+\sqrt{1-e\cos{\theta}}} \right)^3 \sin{\theta} \, d\theta
    \\
    =-\frac{4\pi a^3}{3e} \left(
    \frac{2\sqrt{1-e}+1}{2(\sqrt{1-e}+1)^2}-
        \frac{2\sqrt{1+e}+1}{2(\sqrt{1+e}+1)^2}
\right).
    \end{split}
\end{align}
Also, moment of the shape is
\begin{align}
    \begin{split}
        M=\frac{\pi}{2}\int \limits_{0}^{\pi} \left( \frac{a}{1+\sqrt{1-e\cos{\theta}}} \right)^4 \cos{\theta} \sin{\theta} \, d\theta
        \\
        =\frac{\pi a^4}{e^2}\left(
        \ln{\left(\frac{\sqrt{1+e}+1}{\sqrt{1-e}+1}\right)}+
        \frac{3\sqrt{1+e}+2}{(\sqrt{1+e}+1)^2}-\frac{3\sqrt{1-e}+2}{(\sqrt{1-e}+1)^2}
        \right).
    \end{split}
\end{align}
Furthermore, the average cross sectional area is
\begin{align}
    \begin{split}
        T(2,3)=\frac{\pi}{2}\int \limits_{0}^{\pi} \left( \frac{a}{1+\sqrt{1-e\cos{\theta}}} \right)^2 \sin{\theta} \, d\theta
        \\
        =-\frac{\pi a^2}{e}\left( 
        \ln{\left( \frac{1+\sqrt{1+e}}{1+\sqrt{1-e}}\right)}+\frac{1}{1+\sqrt{1+e}}-\frac{1}{1+\sqrt{1-e}}
        \right).
    \end{split}
\end{align}
Integration in $n=3$, $m=2$ is straight forward for every cases, whereas only properties of critical shape is integrated in other cases.
\subsection{n = 5, m = 1}
The equation of the stationary shape is a quintic equation, so the equation cannot be solved algebraically. However, it is still possible to obtain the hyper cylindrical equation and integrate the volume and the moment of the critical shape. From equation (6.8),
\begin{equation}
    R^2=\sqrt{\frac{1}{\lambda+\mu z}}-z^2.
\end{equation}
The volume and the moment is integrated in hyper cylindrical coordinates as
\begin{align}
    \begin{split}
        V=\int \limits_{z_{-}}^{z_+} V_4 R^4 \, dz
        \\
        =\frac{\pi^2}{2}\left[
    \frac{\ln{(\lambda+\mu z)}}{\mu}+\frac{z^5}{5}-\frac{\sqrt{\lambda+\mu z}\,(12\mu^2z^2-16\lambda \mu z + 32\lambda^2)}{15\mu^3}
    \right]_{z_{-}}^{z_+}    
    \end{split}
\end{align}
and
\begin{align}
    \begin{split}
        M=\int \limits_{z_{-}}^{z_+} V_4 R^4 z\, dz =\frac{\pi^2}{2}
    \biggl[
    \frac{-\lambda \ln{(\lambda+\mu z)}}{\mu^2}+\frac{z^6}{6}
    \\
    -\frac{\sqrt{\lambda+\mu z}\,(120\mu^3z^3-144\lambda \mu^2 z^2 + 192\lambda^2\mu z -384 \lambda^3)}{210\mu^3}
    \biggr]_{z_{-}}^{z_+}.
    \end{split}
\end{align}
where $z_+$ and $z_-$ can be solved since the hyper cylindrical equation has double root at $z_+$. The average volume of the cross section can be obtained through the equation (7.6). $R=0$ gives
\begin{equation}
    5w^4-4w^5=1
\end{equation}
where $w=z\left(\frac{\lambda}{5} \right)^{\frac{1}{4}}$. The equation (7.33) has a trivial solution $w=1$, and solving the rest of the equation through Cardano method gives $z_+=\left(\frac{5}{\lambda} \right)^{\frac{1}{4}}$ and
\newline $z_-=-\left(\frac{5}{\lambda}\right)^{\frac{1}{4}}\frac{1}{12}\left(\sqrt[3]{15(4\sqrt{6}-9)}-\sqrt[3]{15(4\sqrt{6}+9)}-3\right)$.
\subsection{n - m is odd}
For $n-m\geq 3$ and odd, equation of the stationary shape becomes more complicated as we can notice in the equation (6.8). However, it is possible to integrate the volume and the moment of the critical shape for $n-m \leq 5$ in the hyper cylindrical coordinate. To solve for the limits of the integral $z_+$ and $z_-$, set $R=0$, and substitute the equation (7.4) to obtain
\begin{equation}
    1-\lambda z^{n-m}- \frac{\lambda^{\frac{n-m+1}{n-m}}}{(n-m+1)^{\frac{n-m+1}{n-m}}}(n-m) z^{n-m+1}=0.
\end{equation}
Substituting $w=\frac{\lambda^{\frac{1}{n-m}}}{(n-m+1)^{\frac{1}{n-m}}}z$ gives a simplified equation
\begin{align}
    \begin{split}
    1-(n-m+1)w^{n-m}+(n-m)w^{n-m+1}
    \\
    =(w-1)^2\left(
    (n-m)w^{n-m-1}+(n-m-1)w^{n-m-2}+ \cdots + 2w +1
    \right)=0.
    \end{split}
\end{align}
For $n-m \leq 5$, the equation, and the limits of the integration $z_+$ and $z_-$, can be solved algebraically. Integration in hyper cylindrical coordinates give the volume and the moment of the critical shape as
\begin{equation}
    V=\int \limits_{z_{-}}^{z_+} V_{n-1} R^{n-1} \, dz
\end{equation}
and
\begin{equation}
    M=\int \limits_{z_{-}}^{z_+} V_{n-1} R^{n-1} z\, dz.
\end{equation}
The average thickness of the cross section is calculated through the equation (7.6).
\section{Conclusion}
The average $m$-dimensional cross section is maximized at the sphere, when the centre of mass is constrained to be at the origin. However, the average $m$-dimensional cross section can not be maximized in the extended problem where the centre of mass is not at the origin. In this case, the shapes which exhibit stationary values of average cross sections exist.

Spherical coordinates in $n$-dimension, and the volume and the surface area of $n$-dimensional sphere is studied, in order to investigate the shapes in arbitrary dimensions. Average thickness of $m$-dimensional cross section in $n$-dimensional shape is examined. Due to spherical symmetry, the average thickness is a simple integration about $(n-1)$-dimensional solid angle. We first started to obtain the formula for the average thickness with an assumption of the shape being star shaped, and we extended the formula such that it can be applied to any shapes.

The shape is solved through calculus of variation, which actually gives the shape with stationary average thickness rather than maximizing it. The stationary shapes are obtained by two method. First, under the assumption that the shape is star shaped, Euler-Lagrange equation with Lagrange multiplier method is used. Second, we give a certain deformation to the stationary shapes, and apply the principles of calculus of variation. Then, the shapes are solved with the principles of linear algebra.

The axial symmetry of the stationary shape is found. The axial symmetry simplifies the equation of the shape. The equation of the stationary shape turns out to be only dependent on the difference of the dimension of the cross section and the shape. Also, the axial symmetry let us to use the hyper cylindrical coordinates to describe the stationary shapes. The hyper cylindrical coordinates show that the stationary shapes are either the egg shape, the critical shape, or the opened region.

The relationship between the constants of the equation and the properties of the shapes, such as volume, moment, and average thickness, are obtained using integration. Only $n=3$, $m=2$ case is integrated to elementary functions. In other cases, the critical shapes can be integrated into elementary functions. The spherical equations cannot be solved algebraically if $n-m \geq 5$, and integration of the volume and the moment of critical shapes in elementary functions is possible for $n-m \leq 5$.


\acknowledgments

I thank professor Andrew Hodges for supervising the project and useful discussions.


\begin{thebibliography}{99}

\bibitem {A} Blasjo, Viktor, \textit{The isoperimetric problem},
The American Mathematical Monthly 112.6: 526-566, 2005.

\bibitem {M} Blumenson,\textit{A Derivation of n-Dimensional Spherical Coordinates}, The American Mathematical Monthly 67: 63–66, 1960.

\bibitem {B} Osserman, \textit{The isoperimetric inequality}, Bulletin of the American Mathematical Society 84.6: 1182-1238, 1978.

\bibitem {C} Hilbert, David, and Stephan Cohn-Vossen, \textit{Geometry and the Imagination}, No. 87. American Mathematical Soc., 1999.

\bibitem {D} Weisstein, Eric W, \textit{Hypersphere}, URL http://mathworld.wolfram.com/Hypersphere.html, 2014.

\bibitem {E} Huber, Greg, \textit{Gamma function derivation of n-sphere volumes}, The American Mathematical Monthly 89.5: 301-302, 1982.

\bibitem {F} Weisstein, Eric W, \textit{Gamma Function}, URL http://mathworld.wolfram.com/GammaFunction. html, 2012.

\bibitem {G} Weisstein, Eric W, \textit{Gaussian Integral}, URL http://mathworld.wolfram.com/GaussianIntegral. html, 2004.

\bibitem {H} Weisstein, Eric W, \textit{Lagrange Multiplier}, URL http://mathworld.wolfram.com/LagrangeMulti plier.html, 1999.

\bibitem {I} Lang, Serge, \textit{Calculus of several variables}, Springer Science \& Business Media, 2012.

\bibitem {J} Hodges, Andrew, \textit{Calculus of Variation}, University of Oxford, 2016.

\bibitem {K} Arfken, George B., and Hans J. Weber, \textit{Mathematical methods for physicists}, 1999.




\end{thebibliography}
\end{document}